\documentclass[12pt]{amsart}
\usepackage{mathrsfs}
\usepackage{amsfonts}
\usepackage{amssymb}
\usepackage{amsfonts, amscd, amsmath, mathrsfs, amssymb, amsthm, amsxtra, stmaryrd, bbding, epsfig, graphicx, latexsym, url, mathbbol, bbold}
\usepackage[papersize={7.5in,10.5in},textwidth=14.9cm,textheight=21cm,centering]{geometry}
\usepackage{enumerate}

\usepackage[colorlinks=true,citecolor=blue,linkcolor=blue]{hyperref}
\hypersetup{
pdfstartpage=1,
pdfstartview=FitH}

\DeclareFontFamily{U}{mathx}{\hyphenchar\font45}
\DeclareFontShape{U}{mathx}{m}{n}{
      <5> <6> <7> <8> <9> <10>
      <10.95> <12> <14.4> <17.28> <20.74> <24.88>
      mathx10
      }{}
\DeclareSymbolFont{mathx}{U}{mathx}{m}{n}
\DeclareMathAccent{\widecheck}{\mathalpha}{mathx}{"71}

 \usepackage{caption} 
\numberwithin{equation}{section}

\allowdisplaybreaks

\newtheorem{theorem}{Theorem}[section]
\newtheorem{lemma}{Lemma}[section]

\makeatletter
\newcounter{roem}
\renewcommand{\theroem}{\Roman{roem}}

\newcommand{\c@org@eq}{}
\let\c@org@eq\c@equation
\newcommand{\org@theeq}{}
\let\org@theeq\theequation

\newcommand{\setroem}{
\let\c@equation\c@roem
 \let\theequation\theroem}

\newcommand{\setarab}{
\let\c@equation\c@org@eq
\let\theequation\org@theeq}
\makeatother

\newtheorem*{claim*}{Claim}

\theoremstyle{remark}

\newcommand{\ue}{\mathrm{e}}

\newcommand{\Tr}{\mathrm{Tr}}

\newcommand{\kl}{\mathrm{Kl}}

\newcommand{\bh}{\mathbf{h}}

\newcommand{\bx}{\mathbf{x}}
\newcommand{\by}{\mathbf{y}}
\newcommand{\bz}{\mathbf{z}}

\newcommand{\bC}{\mathbf{C}}
\newcommand{\bF}{\mathbf{F}}
\newcommand{\bN}{\mathbf{N}}

\newcommand{\bZ}{\mathbf{Z}}

\newcommand{\cB}{\mathcal{B}}

\newcommand{\cH}{\mathcal{H}}

\newcommand{\cM}{\mathcal{M}}
\newcommand{\cN}{\mathcal{N}}

\newcommand{\cS}{\mathcal{S}}

\newcommand{\cX}{\mathcal{X}}
\newcommand{\cY}{\mathcal{Y}}

\newcommand{\fq}{\mathfrak{q}}

\newcommand{\rF}{\mathrm{F}}

\newcommand{\bdeta}{\boldsymbol\eta}

\newcommand{\bdkappa}{\boldsymbol\kappa}
\newcommand{\bdlambda}{\boldsymbol\lambda}
\newcommand{\bdchi}{\boldsymbol\chi}
\newcommand{\bdrho}{\boldsymbol\rho}

\usepackage{graphicx}
\usepackage{tikz}

\begin{document}

\vglue -2mm

\title{Equidistributions of Jacobi sums}

\author{Ping Xi}

\address{School of Mathematics and Statistics, Xi'an Jiaotong University, Xi'an 710049, P. R. China}

\email{ping.xi@xjtu.edu.cn}



\begin{abstract}
Let $\bF_q$ be a finite field of $q$ elements. We show that the normalized Jacobi sum
$J(\chi,\eta)/\sqrt{q}$, for each fixed non-trivial multiplicative character $\eta$, becomes equidistributed in the unit circle as $q\rightarrow+\infty,$ when $\chi$ runs over all 
non-trivial multiplicative characters different from $\eta^{-1}.$
Previously, the similar equidistribution was obtained by Katz and Zheng by varying both of $\chi$ and $\eta$. On the other hand, we also obtain the equidistribution of $J(\chi,\eta)$ as $(\chi,\eta)$ runs over $\cX\times\cY\subseteq(\widehat{\rF^*})^2$, as long as $|\cX|>q^{\frac{1}{2}+\varepsilon}$ and $|\cY|>q^\varepsilon$ for any $\varepsilon>0$. This updates a recent work of Lu, Zheng and Zheng, who require $|\cX||\cY|>q\log^2q.$

The main ingredient is the estimate for hypergeometric sums due to Katz.
\end{abstract}
\vglue -3mm
\maketitle

\setcounter{tocdepth}{1}

\section{Introduction}
Let $\rF=\bF_q$ be a finite field of characteristic $p$ with $q$ elements. Given two multiplicative characters $\chi$ and $\eta$, the Jacobi sum can be defined by
\begin{align*}
J(\chi,\eta)=\sum_{x\in\rF^*}\chi(x)\eta(1-x).
\end{align*}
It is known that $|J(\chi,\eta)|=\sqrt{q}$ if $\chi,\eta$ and $\chi\eta$ are all non-trivial. 
The main concern is to study how the quotient $J(\chi,\eta)/\sqrt{q}$, for each fixed non-trivial character $\eta$, is distributed in the unit circle as long as $\chi$ runs over all characters in $G(\eta):=\widehat{\rF^*}\setminus\{\mathbf{1},\eta^{-1}\}.$ The relevant problem was considered by Katz for Gauss sums
 \begin{align*}
\tau(\chi)=\sum_{x\in\rF^*}\chi(x)\psi(x),
\end{align*}
where $\psi$ is an additive character given by $\psi(x)=\exp(2\pi i\Tr_{\bF_q/\bF_p}(x)/p).$
More precisely, Katz \cite{Ka80}, using Deligne's estimate for hyper-Kloosterman sums, proved that $\psi(\chi)/\sqrt{q}$ is equidistributed on the unit circle as $\chi$ runs over $\widehat{\rF^*}.$ The case of Jacobi sums was considered by Katz and Zheng \cite{KZ96}, who proved that $J(\chi,\eta)/\sqrt{q}$ is equidistributed on the unit circle as both of $\chi$ and $\eta$ run over $\widehat{\rF^*}$ with $\chi,\eta,\chi\eta$ non-trivial. The first aim of this paper is to show that the equidistribution is also valid if one character is fixed.

\begin{theorem}\label{thm:equidistribution}
Let $\eta$ be a non-trivial multiplicative character in $\widehat{\rF^*}$. For any given interval $I\subseteq~]-\pi,\pi],$ we have
\begin{align*}
\frac{|\{\chi\in G(\eta):\arg J(\chi,\eta)\in I\}|}{q-3}=\frac{|I|}{2\pi}+O(q^{-\frac{1}{4}}).
\end{align*}
\end{theorem}

Moreover, we have the following joint distributions for Jacobi sums.

\begin{theorem}\label{thm:jointequidistribution}
Let $s\geqslant2$ be a positive integer. Let $\bdeta=(\eta_1,\eta_2,\cdots,\eta_s)\in (\widehat{\rF^*}\setminus\{\mathbf{1}\})^s$ be a tuple of pairwise distinct multiplicative characters. For any given box $\cB\subseteq~]-\pi,\pi]^s,$ we have
\begin{align*}
\frac{|\{\chi\in G(\bdeta):\{\arg J(\chi,\eta_i)\}_{1\leqslant i\leqslant s}\in\cB\}|}{q-s-2}=\frac{|\cB|}{(2\pi)^s}+O(3^sq^{-\frac{1}{2(s+1)}}).
\end{align*}
where $G(\bdeta)=\bigcap_{1\leqslant i\leqslant s}G(\eta_i)$ and the implied constant is absolute.
\end{theorem}

Let $r,s$ be two fixed positive integers. We now define the mixed moment
\begin{align}\label{eq:definition-M}
\cM_{\bdkappa,\bdlambda}(\bdeta,\bdrho)
&=\sum_{\chi\in G(\bdeta,\bdrho)}\prod_{1\leqslant i\leqslant r}J(\chi,\eta_i)^{\kappa_i}\prod_{1\leqslant j\leqslant s}\overline{J(\chi,\rho_j)}^{\lambda_j},
\end{align}
where
\begin{itemize}
\item $\bdkappa=(\kappa_1,\kappa_2,\cdots,\kappa_r)\in\bN^r$, $\bdlambda=(\lambda_1,\lambda_2,\cdots,\lambda_s)\in\bN^s$ are two non-zero vectors;

\item $\bdeta=(\eta_1,\eta_2,\cdots,\eta_r)\in(\widehat{\rF^*})^r$ and $\bdrho=(\rho_1,\rho_2,\cdots,\rho_s)\in(\widehat{\rF^*})^s$ are two tuples of multiplicative characters with 
$\eta_i\neq\rho_j$ for all $1\leqslant i\leqslant r$ and $1\leqslant j\leqslant s,$ in which case we say $\bdeta$ and $\bdrho$ are {\it disjoint};

\item $G(\bdeta,\bdrho)=G(\bdeta)\cap G(\bdrho).$
\end{itemize}
If $s=0$, the characters $\bdrho$ will not appear, and we denote the relevant moment by $\cM_{\bdkappa}(\bdeta)$.

For $\chi\in G(\bdeta,\bdrho)$, each Jacobi sum in \eqref{eq:definition-M} is of modulus $\sqrt{q}$, and it follows trivially that
\begin{align}\label{eq:trivialbound}
|\cM_{\bdkappa,\bdlambda}(\bdeta,\bdrho)|\leqslant q^{\frac{\|\bdkappa\|_1+\|\bdlambda\|_1}{2}+1}.
\end{align}

To obtain the above equidistributions, we have to beat the trivial bound in \eqref{eq:trivialbound}, as presented by the following theorem.

\begin{theorem}\label{thm:moments}
Let $\bdkappa\in\bZ_+^r,\bdlambda\in\bZ_+^s$ with $\kappa=\|\bdkappa\|_1$ and $\lambda=\|\bdlambda\|_1.$ Assume the two tuples of non-trivial multiplicative characters $\bdeta$ and $\bdrho$ are disjoint. Then we have
\begin{align*}
|\cM_{\bdkappa,\bdlambda}(\bdeta,\bdrho)|\leqslant 3(\kappa+\lambda)q^{\frac{\kappa+\lambda+1}{2}}.
\end{align*}
\end{theorem}

Let $\cX,\cY\subseteq\widehat{\rF^*}$ be two sets of non-trivial multiplicative characters. For each positive integer $\kappa$, we define
the bilinear moment
\begin{align*}
\cN_\kappa(\cX,\cY):=\mathop{\sum_{\chi\in\cX}\sum_{\eta\in \cY}}_{\chi\eta\neq1}J(\chi,\eta)^\kappa.\end{align*}

\begin{theorem}\label{thm:bilinearmoment}
With the above notation, we have
\begin{align*}
|\cN_\kappa(\cX,\cY)|
&\leqslant s|\cX|^{1-\frac{1}{2s}}(|\cY|^sq^{s\kappa+1}+6s\kappa|\cY|^{2s}q^{s\kappa+\frac{1}{2}})^{\frac{1}{2s}}\end{align*}
for all positive integers $\kappa$ and $s.$
\end{theorem}

From Theorem \ref{thm:bilinearmoment}, we may conclude the following bilinear equidistribution of Jacobi sums.

\begin{theorem}\label{thm:bilinearequidistribution}
Let $\cX,\cY\subseteq\widehat{F^*}$ be two sets of multiplicative characters in $\widehat{\rF^*}$ with $|\cX|\geqslant\sqrt{q}$ and $|\cY|\geqslant2.$ Let $s\geqslant1$ be a positive integer. For any given interval $I\subseteq~]-\pi,\pi],$ we have
\begin{align*}
\frac{|\{(\chi,\eta)\in\cX\times\cY:\arg J(\chi,\eta)\in I,\chi\eta\neq\mathbf{1}\}|}{|\cX||\cY|}=\frac{|I|}{2\pi}+O\Big\{s\Big(\frac{\sqrt{q}}{|\cX|}\Big)^{\frac{1}{2s+1}}+s|\cY|^{-\frac{1}{2}}\Big(\frac{q\log q}{|\cX|}\Big)^{\frac{1}{2s}}\Big\},
\end{align*}
where the implied constant is absolute. 
\end{theorem}

In many problems of analytic number theory, one has to control the bilinear form
\[\sum_{m\leqslant M}\sum_{n\leqslant N}\alpha_m\beta_n W_\fq(mn)\]
non-trivially, where $(\alpha_m)_{m\leqslant M}$ and $(\beta_n)_{n\leqslant N}$ are arbitrary complex coefficients, and $W_\fq:\bZ/\fq\bZ\rightarrow\bC$ is the arithmetic object which we are interested in. Due to the structure of intervals and by completion, one can usually work non-trivially as long as 
$M>\fq^{\frac{1}{2}+\varepsilon}$ and $N>\fq^{\varepsilon}$, provided that $W_\fq$ allows certain oscillations against additive characters.
The situation in Theorem \ref{thm:bilinearmoment} is quite different; $\cX,\cY$ have no special structures such as ``intervals", and the completion method is not available in such case. Fortunately, by taking $s$ sufficiently large in Theorems \ref{thm:bilinearmoment} and \ref{thm:bilinearequidistribution}, we find that the estimate in Theorem \ref{thm:bilinearmoment} is non-trivial as long as $|\cX|>q^{\frac{1}{2}+\varepsilon}$ and $|\cY|>q^{\varepsilon}$ for any given $\varepsilon>0,$
in which case Theorem \ref{thm:bilinearequidistribution} gives the desired equidistribution of $J(\chi,\eta)$ as $\chi$ and $\eta$ run over $\cX$ and $\cY$, respectively.
In fact, equidistribution in Theorem \ref{thm:bilinearequidistribution} was asked by Shparlinski \cite{Sh09}, and the first attack is due to Lu, Zheng and Zheng \cite{LZZ18}, who require that $|\cX||\cY|>q\log^2q.$ 

Before closing this section, we would like to mention that it is interesting to extend all the above results to the high-dimensional Jacobi sum
\[J(\chi_1,\chi_2,\cdots,\chi_k)=\mathop{\sum\sum\cdots\sum}_{\substack{x_1,x_2,\cdots,x_k\in\rF^*\\x_1+x_2+\cdots+x_k=1}}\chi_1(x_1)\chi_2(x_2)\cdots \chi_k(x_k)\]
as already considered in \cite{LZZ18}, where $\chi_i\in\widehat{\rF^*}$ for $1\leqslant i\leqslant k.$

\subsection*{Acknowledgements} 
The work was inspired by a conversation with Professor Wen-Ching Winnie Li during my visit to Institute of Mathematical Research, The University of Hong Kong (HKU) in the summer of 2018. I would like to thank Professor Li for her inspiration, to Dr. Yuk-Kam Lau for his kind invitation, and to HKU for the support and hospitality.
The work is supported in part by NSFC (No.11601413) and NSBRP (No. 2017JQ1016) of Shaanxi Province.

\smallskip

\section{Preparations}

The following lemma allows us to express Jacobi sums in terms of Gauss sums.
\begin{lemma}\label{lm:Jacobi-Gauss}
Suppose $\chi,\eta$ and $\chi\eta$ are all non-trivial multiplicative characters in $\widehat{\rF^*}.$
Then we have 
\begin{align*}
J(\chi,\eta)=\frac{\tau(\chi)\tau(\eta)}{\tau(\chi\eta)}.
\end{align*}
\end{lemma}

The tools to finish our estimates are hypergeometric sums introduced by Katz (see \cite[Chapter 8]{Ka90}). Let $m,n$ be two non-negative integers. Suppose $\boldsymbol{\chi}=(\chi_i)_{1\leqslant i\leqslant m}$ and $\boldsymbol{\eta}=(\eta_j)_{1\leqslant j\leqslant n}$ are two tuples of multiplicative characters in $\widehat{\rF^*}$. Katz introduced the following hypergeometric sum
\begin{align*}
H(t;\boldsymbol\chi,\boldsymbol\eta,q)
:=\frac{(-1)^{m+n-1}}{q^{(m+n-1)/2}}\mathop{\sum\sum}_{\substack{\bx\in(\rF^*)^m,\by\in(\rF^*)^n\\ N(\bx)=tN(\by)}}\boldsymbol\chi(\bx)
\overline{\boldsymbol\eta(\by)}\psi(T(\bx)-T(\by))
\end{align*}
for $t\in\rF^*,$ where, for $\bx=(x_1,x_2,\cdots,x_m)\in(\rF^*)^m$, 
\begin{align*}
\boldsymbol\chi(\bx)=\prod_{1\leqslant i\leqslant m}\chi_i(x_i),
\end{align*}
\begin{align*}
T(\bx)=x_1+x_2+\cdots+x_m,\ \ \ N(\bx)=x_1x_2\cdots x_m,
\end{align*}
and the notation with $\by$ can be defined in the same way.

As one may see, $H(t;\boldsymbol\chi,\boldsymbol\eta,q)$ extends the classical (hyper-) Kloosterman sums. In particular, if $n=0$ and $\chi_i=\mathbf{1}$ for each $1\leqslant i\leqslant m,$ the hypergeometric sum gives the hyper-Kloosterman sum
\begin{align*}
\kl_m(t,q)=\frac{(-1)^{m-1}}{q^{(m-1)/2}}\sum_{\substack{\bx\in(\rF^*)^m\\ N(\bx)=t}}\psi(T(\bx)),
\end{align*}
which was already proved by Deligne \cite{De80} that $t\mapsto\kl_m(t,q)$ is the Frobenius trace function of a certain $\ell$-adic sheaf ($\ell\neq p$), which is geometrically irreducible of rank $m$, lisse on $\mathbf{G}_{m,\bF_q}$ and pointwise pure of weight 0.

In general, Katz \cite[Theorem 8.4.2]{Ka90} proved the following assertion, which is quite crucial in proving equidistributions as above.

\begin{lemma}\label{lm:hypergeometricsum}
With the above notation,
if $\bdchi$ and $\bdeta$ are disjoint, then for any $\ell\neq p,$ there exists a geometrically irreducible $\ell$-adic middle-extension sheaf $\cH(\boldsymbol\chi,\boldsymbol\eta)$ on $\mathbf{A}_{\bF_q}^1$ with trace function given by $t\mapsto H(t;\bdchi,\bdeta,q),$ such that it is
\begin{itemize}
\item pointwise pure of weight $0$ and of rank $\max\{m,n\};$
\item lisse on $\mathbf{G}_{m,\bF_q}$, if $m\neq n;$
\item lisse on $\mathbf{G}_{m,\bF_q}-\{1\}$ and of rank $m,$ if $m=n.$
\end{itemize}
\end{lemma}

\smallskip

\section{Concluding equidistributions}
We first give the proof of Theorem \ref{thm:moments}.
The proof of Theorems \ref{thm:equidistribution} and \ref{thm:jointequidistribution} will be the same, and we only give the details for the latter one. 
\subsection{Proof of Theorem \ref{thm:moments}} We write $\kappa=\|\bdkappa\|_1$ and $\lambda=\|\bdlambda\|_1.$ Without loss of generality, we assume that $\kappa\geqslant\lambda.$
From Lemma \ref{lm:Jacobi-Gauss} it follows that
\begin{align*}
\cM_{\bdkappa,\bdlambda}(\bdeta,\bdrho)
&=q^{-\kappa-\lambda}
\sum_{\chi\in G(\bdeta,\bdrho)}\tau(\chi)^\kappa\overline{\tau(\chi)}^\lambda\prod_{1\leqslant i\leqslant r}\tau(\eta_i)^{\kappa_i}\overline{\tau(\chi\eta_i)}^{\kappa_i}\cdot \prod_{1\leqslant j\leqslant s}\overline{\tau(\rho_j)}^{\lambda_j}\tau(\chi\rho_j)^{\lambda_j}\\
&=\vartheta(\bdeta,\bdrho)\cdot q^{\frac{\lambda-\kappa}{2}}
\sum_{\chi\in G(\bdeta,\bdrho)}\tau(\chi)^{\kappa-\lambda}\prod_{1\leqslant i\leqslant r}\overline{\tau(\chi\eta_i)}^{\kappa_i}\cdot \prod_{1\leqslant j\leqslant s}\tau(\chi\rho_j)^{\lambda_j},
\end{align*}
where $\vartheta(\bdeta,\bdrho)$ is a constant depending only on $\bdeta,\bdrho$ with $|\vartheta(\bdeta,\bdrho)|=1.$

Note that 
\begin{align*}
\Bigg|\sum_{\chi\in \widehat{\rF^*}\setminus G(\bdeta,\bdrho)}\tau(\chi)^{\kappa-\lambda}\prod_{1\leqslant i\leqslant r}\overline{\tau(\chi\eta_i)}^{\kappa_i}\cdot \prod_{1\leqslant j\leqslant s}\tau(\chi\rho_j)^{\lambda_j}\Bigg|\leqslant \sum_{1\leqslant i\leqslant r}q^{\kappa-\frac{\kappa_i}{2}}+\sum_{1\leqslant j\leqslant s}q^{\lambda-\frac{\lambda_j}{2}}.
\end{align*}
Therefore, we may write
\begin{align}\label{eq:M-M*}
|\cM_{\bdkappa,\bdlambda}(\bdeta,\bdrho)-\cM^*_{\bdkappa,\bdlambda}(\bdeta,\bdrho)|\leqslant q^{\frac{\lambda-\kappa}{2}}\Big(\sum_{1\leqslant i\leqslant r}q^{\kappa-\frac{\kappa_i}{2}}+\sum_{1\leqslant j\leqslant s}q^{\lambda-\frac{\lambda_j}{2}}\Big),
\end{align}
where
\begin{align*}
\cM^*_{\bdkappa,\bdlambda}(\bdeta,\bdrho)
&=\vartheta(\bdeta,\bdrho)\cdot q^{\frac{\lambda-\kappa}{2}}
\sum_{\chi\in \widehat{\rF^*}}\tau(\chi)^{\kappa-\lambda}\prod_{1\leqslant i\leqslant r}\overline{\tau(\chi\eta_i)}^{\kappa_i}\cdot \prod_{1\leqslant j\leqslant s}\tau(\chi\rho_j)^{\lambda_j}.
\end{align*}
Opening Gauss sums and
by virtue of the orthogonality of additive characters, we derive that
\begin{align*}
&\ \ \ \ \cM^*_{\bdkappa,\bdlambda}(\bdeta,\bdrho)\\
&=\vartheta(\bdeta,\bdrho)\cdot q^{\frac{\lambda-\kappa}{2}}
\mathop{\sum\cdots\sum}_{\substack{\bx\in(\rF^*)^{\kappa-\lambda}\\ \by_1\in(\rF^*)^{\kappa_1},\cdots,\by_r\in(\rF^*)^{\kappa_r}\\
\bz_1\in(\rF^*)^{\lambda_1},\cdots,\bz_s\in(\rF^*)^{\lambda_s}}}\psi\Big(T(\bx)-\sum_{1\leqslant i\leqslant r}T(\by_i)+\sum_{1\leqslant j\leqslant s}T(\bz_j)\Big)\\
&\ \ \ \ \times\sum_{\chi\in \widehat{\rF^*}}\chi(N(\bx))\cdot\prod_{1\leqslant i\leqslant r}\overline{(\chi\eta_i)(N(\by_i))}\cdot\prod_{1\leqslant j\leqslant s}(\chi\rho_j)(N(\bz_j))\\
&=\vartheta(\bdeta,\bdrho)\cdot (q-1)q^{\frac{\lambda-\kappa}{2}}
\mathop{\sum\cdots\sum}_{\substack{\bx\in(\rF^*)^{\kappa-\lambda}\\ \by_1\in(\rF^*)^{\kappa_1},\cdots,\by_r\in(\rF^*)^{\kappa_r}\\
\bz_1\in(\rF^*)^{\lambda_1},\cdots,\bz_s\in(\rF^*)^{\lambda_s}\\
N(\bx)\prod_{1\leqslant j\leqslant s}N(\bz_j)=\prod_{1\leqslant i\leqslant r}N(y_i)}}\prod_{1\leqslant i\leqslant r}\overline{\eta_i(N(\by_i))}\cdot\prod_{1\leqslant j\leqslant s}\rho_j(N(\bz_j))\\
&\ \ \ \ \times\psi\Big(T(\bx)-\sum_{1\leqslant i\leqslant r}T(\by_i)+\sum_{1\leqslant j\leqslant s}T(\bz_j)\Big).
\end{align*}
Note that the multiple sum over $\bx,\by_1,\cdots,\by_r,\bz_1,\cdots,\bz_s$ is exactly
\begin{align*}
-q^{\kappa-\frac{1}{2}}H(1; (\underbrace{\mathbf{1},\cdots,\mathbf{1}}_{\kappa-\lambda\text{ copies}}, \underbrace{\rho_1,\cdots,\rho_1}_{\lambda_1\text{ copies}},\cdots,\underbrace{\rho_s,\cdots,\rho_s}_{\lambda_s\text{ copies}}),(\underbrace{\eta_1,\cdots,\eta_1}_{\kappa_1\text{ copies}},\cdots,\underbrace{\eta_r,\cdots,\eta_r}_{\kappa_r\text{ copies}}),q).
\end{align*}
By Lemma \ref{lm:hypergeometricsum}, the above hypergeometric sum is bounded by $\kappa$ in absolute values.
Hence
\begin{align*}
|\cM^*_{\bdkappa,\bdlambda}(\bdeta,\bdrho)|
&\leqslant \kappa q^{\frac{\kappa+\lambda+1}{2}},
\end{align*}
from which and \eqref{eq:M-M*}, we have
\begin{align*}
|\cM_{\bdkappa,\bdlambda}(\bdeta,\bdrho)|
&\leqslant \kappa q^{\frac{\kappa+\lambda+1}{2}}+q^{\frac{\lambda-\kappa}{2}}\Big(\sum_{1\leqslant i\leqslant r}q^{\kappa-\frac{\kappa_i}{2}}+\sum_{1\leqslant j\leqslant s}q^{\lambda-\frac{\lambda_j}{2}}\Big)\\
&\leqslant 3(\kappa+\lambda) q^{\frac{\kappa+\lambda+1}{2}}
\end{align*}
as stated.

\subsection{Proof of Theorem \ref{thm:jointequidistribution}}
To conclude Theorem \ref{thm:jointequidistribution}, we would like to introduce the notation for {\it discrepancy} of dimension $r\geqslant1$. Let $\cS=\{\bx_1,\bx_2,\cdots,\bx_N\}$ with $\bx_j\in[0,1]^r$ for each $1\leqslant j\leqslant N.$ We define the discrepancy of $\cS$ by
\begin{align*}
D_{N,r}(\cS):=\sup_{B}\Bigg|\frac{1}{N}|\{\bx_n\in B:1\leqslant n\leqslant N\}|-\mathrm{meas}(B)\Bigg|,
\end{align*}
where the supremum is taken over all boxes $B\subseteq[0,1]^r.$
The famous The Erd\H{o}s--Tur\'an--Koksma inequality, as a high-dimensional generalization of the classical Erd\H{o}s--Tur\'an inequality, gives a quantitative form of Weyl's criterion for equidistributions and asserts that, for any $H\geqslant1$,
\begin{align*}
D_{N,r}(\cS)\leqslant 3^r\Big(\frac{2}{H+1}+\sum_{\mathbf{0}\neq\bh\in[0,H]^r\cap\bZ^r}\frac{1}{\beta(\bh)}\Big|\frac{1}{N}\sum_{n\leqslant N}\ue(\langle\bx_n,\bh\rangle)\Big|\Big),
\end{align*}
where $\beta(\bh):=\max_{1\leqslant j\leqslant r}\{1,|h_j|\}$.

Note that 
\[\frac{J(\chi,\eta)}{|J(\chi,\eta)|}=\frac{J(\chi,\eta)}{\sqrt{q}}=\ue^{i\arg J(\chi,\eta)}\]
 for $\chi,\eta,\chi\eta\neq\mathbf{1}.$ Denote by $N_1:=N_1(\bdeta,\cB)$ the number of multiplicative characters counted in Theorem \ref{thm:jointequidistribution}.
Hence, for any $H_1\geqslant3$, we have
\begin{align*}
\frac{N_1}{q-s-2}-\frac{|\cB|}{(2\pi)^s}
&\ll3^s\Bigg(\frac{1}{H_1}+\frac{1}{q}\sum_{\mathbf{0}\neq\bdkappa\in[0,H_1]^s\cap\bZ^s}\frac{1}{\beta(\bdkappa)}\Bigg|\sum_{\chi\in G(\bdeta)}\prod_{1\leqslant j\leqslant s}\ue^{i\kappa_j\arg J(\chi,\eta_j)}\Bigg|\Bigg)\\
&\ll3^s\Big(\frac{1}{H_1}+\frac{1}{q}\sum_{\mathbf{0}\neq\bdkappa\in[0,H_1]^s\cap\bZ^s}\frac{|\cM_{\bdkappa}(\bdeta)|}{\beta(\bdkappa)}q^{-\frac{\|\bdkappa\|_1}{2}}\Big),
\end{align*}
where the implied constants are absolute.
By virtue of Theorem \ref{thm:moments}, we find
\begin{align*}
\frac{N_1}{q-s-2}-\frac{|\cB|}{(2\pi)^s}
&\ll  3^s\Big(\frac{1}{H_1}+\frac{H_1^s}{\sqrt{q}}\Big).
\end{align*}
Then Theorem \ref{thm:jointequidistribution} follows by taking $H_1=q^{\frac{1}{2(s+1)}}$.

\section{Bilinear moments and equidistributions}
\subsection{Proof of Theorem \ref{thm:bilinearmoment}}

By H\"older's inequality, for $s\geqslant1$ we have
\begin{align}\label{eq:Holder}
|\cN_\kappa(\cX,\cY)|^{2s}\leqslant X^{2s-1}\sum_{\chi\in\widehat{\rF^*}}\Bigg|\sum_{\eta\in \cY\setminus\{\chi^{-1}\}}J(\chi,\eta)^\kappa\Bigg|^{2s}=X^{2s-1}\mathop{\sum\sum}_{\bdeta,\bdrho\in\cY^s}\cS(\bdeta,\bdrho),\end{align}
where, for $\bdeta=(\eta_1,\eta_2,\cdots,\eta_s)\in\cY^s,\bdrho=(\rho_1,\rho_2,\cdots,\rho_s)\in\cY^s$,
\begin{align*}
\cS(\bdeta,\bdrho)&=\sum_{\chi\in G(\bdeta,\bdrho)}\prod_{1\leqslant i\leqslant s}J(\chi,\eta_i)^{\kappa}\overline{J(\chi,\rho_i)}^{\kappa}.\end{align*}
It is clear that
\begin{align*}
\cS(\bdeta,\bdrho)&=\cM_{\bdkappa,\bdkappa}(\bdeta,\bdrho),\end{align*}
where $\bdkappa=(\kappa,\kappa,\cdots,\kappa)$.

It suffices to bound $\cS(\bdeta,\bdrho)$ from above. We consider two extreme cases: 

(i) $\bdeta$ coincides with $\bdrho$;

(ii) all coordinates of $\bdeta$ and $\bdrho$ are distinct.

\noindent In the first case, we employ the trivial bound, i.e.,
\begin{align*}
|\cS(\bdeta,\bdrho)|\leqslant q^{s\kappa+1}.\end{align*}
In the second case, we are in a good position to apply Theorem \ref{thm:moments}, getting
\begin{align*}
|\cS(\bdeta,\bdrho)|\leqslant 6s\kappa q^{s\kappa+\frac{1}{2}}.\end{align*}
Note that the above two estimates also hold correspondingly if the coordinates of 
$\bdeta$ and $\bdrho$ are permuted separately. Collecting all possibilities, we have
\begin{align*}
\mathop{\sum\sum}_{\bdeta,\bdrho\in\cY^s}|\cS(\bdeta,\bdrho)|
&\leqslant s^{2s}(Y^{s}q^{s\kappa+1}+6s\kappa Y^{2s}q^{s\kappa+\frac{1}{2}}),\end{align*}
from which and \eqref{eq:Holder} we obtain
\begin{align*}
|\cN_\kappa(\cX,\cY)|
&\leqslant sX^{1-\frac{1}{2s}}(Y^sq^{s\kappa+1}+6s\kappa Y^{2s}q^{s\kappa+\frac{1}{2}})^{\frac{1}{2s}}\end{align*}
as stated.
\subsection{Proof of Theorem \ref{thm:bilinearequidistribution}}
Denote by $N_2:=N_2(\cX,\cY;I)$ the number of pairs of multiplicative characters $(\chi,\eta)$ counted in Theorem \ref{thm:bilinearequidistribution}.
Following the argument to the proof of Theorem \ref{thm:jointequidistribution}, we find, for any $H_2\geqslant3$, that
\begin{align*}
\frac{N_2}{|\cX||\cY|}-\frac{|I|}{2\pi}
&\ll\frac{1}{H_2}+\frac{1}{|\cX||\cY|}\sum_{1\leqslant \kappa\leqslant H_2}\frac{1}{\kappa}\Bigg|\mathop{\sum_{\chi\in \cX}\sum_{\eta\in\cY}}_{\chi\eta\neq1}\ue^{i\kappa\arg J(\chi,\eta)}\Bigg|\\
&\ll\frac{1}{H_2}+\frac{1}{|\cX||\cY|}\sum_{1\leqslant \kappa\leqslant H_2}\frac{|\cN_{\kappa}(\cX,\cY)|}{\kappa}q^{-\frac{\kappa}{2}}.
\end{align*}
By virtue of Theorem \ref{thm:bilinearmoment}, we find
\begin{align*}
\frac{N_2}{|\cX||\cY|}-\frac{|I|}{2\pi}
&\ll\frac{1}{H_2}+s|\cX|^{-\frac{1}{2s}}\{|\cY|^{-\frac{1}{2}}(q\log H_2)^{\frac{1}{2s}}+H_2^{\frac{1}{2s}}q^{\frac{1}{4s}}\}.
\end{align*}
Then Theorem \ref{thm:bilinearequidistribution} follows by taking
\begin{align*}
H_2=2018+(|\cX|/\sqrt{q})^{\frac{1}{2s+1}}.
\end{align*}

\smallskip

\bibliographystyle{plain}

\bigskip

\end{document}